\documentclass[12pt,a4]{amsart}
\usepackage{amsmath, amsthm, amscd, amsfonts, amssymb, graphicx, color}
\usepackage[bookmarksnumbered, colorlinks, plainpages]{hyperref}
 \usepackage{mathrsfs}

\makeatletter \oddsidemargin.9375in \evensidemargin \oddsidemargin
\def\authorsaddresses#1{\dedicatory{#1}}
\newtheorem{theorem}{Theorem}[section]

\theoremstyle{definition}

\newtheorem{example}[theorem]{Example}
\theoremstyle{remark}

\numberwithin{equation}{section}

\begin{document}
\setcounter{page}{1}


\title[Assessing the Rel based on percolation ] {Assessing the reliability polynomial based on percolation theory}

\author[Sajadi]{Sajadi, Farkhondeh A.}

\authorsaddresses{ Department of Statistics, University
of Isfahan, Isfahan 81744,Iran;\\
f.sajadi@sci.ui.ac.ir}


\begin{abstract}
In this paper, we study the robustness of network topologies. We use the concept of percolation as “measuring” tool to assess the reliability polynomial of those systems which can be modeled as a general inhomogeneous random graph as well as scale-free random graph.
\end{abstract}

\maketitle

\section{Introduction}
\label{S:1}
The robustness is one of the structural properties of a complex system which measures its ability of continuing perform well, subject to failures or attacks. It is needed to quantify the measure of robustness in order to decide whether or not, a given system is robust.
The most common measure of robustness of a network to random failures
of components is all-terminal reliability polynomial, the probability that there is an operating communications link between any two components in the system. 
The reliability of a system of interacting agents,  can be determined by analyzing the reliability of the underlying graph. A graph is a pair of two sets $G=(V,E)$ where $V$ is a set of nodes/vertices and $E$ is a set of edges/links that connect two elements of $V$. The elements of a real system and the interactions patterns between them are represented by nodes and edges in a graph, respectability.
Suppose we have a graph $G$ for which the nodes are always operational but for which each edge $e \in E$ is independently operational with probability $p_e \in (0.1)$.  The (all terminal) reliability of $G$, denoted by $Rel(G,p_e)$ is therefore defined to be the probability that the graph is connected when each edge is (independently of the others) present with probability $p_e$. In other words the probability that the operational edges form a spanning connected subgraph of $G$.
Let set of operational edges $S \subseteq E$ be the state of network, i.e., the network is in state $S$ when all edges of $S$ are  operational and all edges $E-S$ are failing. Let $\mathscr{O}$ be the set of all operational states. Then $Rel(G,p_e)$ is equal to  
\begin{equation}
\label{RP}
Rel(G,p_e)=\sum_{S \in \mathscr{O}} \prod_{e \in S} p_{e} \prod_{e \notin S} (1-p_{e})\, .
\end{equation} 
By knowing all states, the reliability is easily but not efficiently computed.
The most amenable case is that in which all of the edge probabilities are identical, say $p$. Therefore under the condition that for all $e \in E, p_e=p$ we  have simpler version of formula(\ref{RP}), which is called F-form of the reliability polynomial, as follows
\begin{equation}
\label{SRP}
Rel(G,p)=\sum_{i=0}^{|E|} F_i (1-p)^{i} p^{|E|-i} 
\end{equation} 
where $F_i$ is the number of connected spanning subgraphs of size $|E|-i$. 
There are different forms of presentation for the reliability polynomial \cite{colb1993}. 
Also there are methods to examine the possible states of the network, for example the factorization method \cite{ch}.
However, computing the reliability polynomial, $Rel(G,p_e)$ is not easy in general. 
Even the time required to calculate $R(G, p)$ for an arbitrary connected graph $G$ grows exponentially with the size of the network, \cite{colb1987},\cite{Pro1983}.

Because of this difficulty, network designers have often relied upon upper and lower bounds on $Rel(G,p)$, \cite{Lomon1971},\cite{Lomon1972},\cite{Kru1963},\cite{Stan1975},\cite{ball1982},\cite{ball1995},\cite{brown1996}. It would be of interest if one can obtain exact solutions for the reliability polynomial  for some families of graphs $G$.
In \cite{daq}, the authors built a framework  based on percolation theory to calculate  the network reliability. They considered homogeneous random graph models with Poisson degree distribution. But in the current work, the method which we used to assess the reliability polynomial can be apply for more general random graph with arbitrary degree distribution. More precisely,    
in this paper, we use the percolation concept as “measuring” tool to assess the reliability polynomial of systems which can be modeled as a general inhomogeneous random graph and scale-free random graph.

The rest of the paper is organized as follows. In Section~\ref{S:2} we describe the inverse percolation process. In Section~\ref{S:3} we explain how to obtain the threshold at which the network loses its connectivity. Section~\ref{S:4} presents the assessment of network reliability. Section~\ref{S:5} and \ref{S:6}, will illustrate our method for two special case of real network which can be modeled as inhomogeneous and scale free random graph, respectively. 

\section{ Inverse percolation process}
\label{S:2}
Random graphs have been used extensively as models for various types of real world networks. It provides techniques to analyze  structure in a system of interacting agents, represented as a network. In mathematics, networks are often referred to as graphs. A random graph is a graph that is sampled according to some probability distribution over  a collection of graphs. 
\textit{Random} graphs are used to prove deterministic properties of the graphs.
The most basic property of a graph is that of being connected. 
It has been shown that, the limiting probability that a random graph possesses connectedness, jumps from $0$ to $1$ (or vice versa) very rapidly. 
There exists a threshold when a transition occurs from not being connected to being connected.
For example for the binomial random graph $\mathcal{G}_{n,p}$, where $n$ is
the number of vertices of the graph and $p \in (0,1)$ is the probability that the edge $(u,v)$ is present, we say that a phase transition occurs if there exists a function $p(n)$ such that for $ p_1(n) << p(n)$, $\mathcal{G}_{n,p_1(n)}$ almost surely does not have the property of being connected, but for $ p(n) << p_2(n)$, $\mathcal{G}_{n,p_2(n)}$ almost surely has the property of being connected. Here $p(n)$ is the threshold.
 The most important transitions is the emergence of a giant component, a connected component of size $\Theta(n)$ where $n$ is the size of $V$, ($a_n=\Theta(b_n)$, if there exist constants $C,c$ and $n_0$, such that, $c b_n \leq a_n \leq C b_n$, for $n\geq n_0$).
The macroscopic behavior of networks, when faced with random  removal of nodes or edges, is
 characterized  in terms of an inverse percolation process in a random graph. Percolation theory characterizes this property of random graphs, especially infinitely large ones. It studies the behavior of the operational giant component, ensuring the global well functioning of the network. Once the random removal of nodes or edges is done (random failure process), the size of giant component, is expected to decrease.
 Using percolation theory, one can predicts the presence of a threshold, above which the network percolates, i.e., it has a giant component and below it the giant component disappears, i.e. the network is fragmented in many disconnected components of very small sizes. For a highly-connected graph sequence, the giant component is unique when the  limit superior of the size of second-largest component by  $n$ is $0$ for large $n$, \cite{Janson}.

\section{ The percolation threshold}
\label{S:3}
Let $p_{e}$ be the probability that edge $e$ is operational and $q_{v}$ be the probability that the node $v \in V$ is operational. By being operational we mean that node/edge  hasn't failed or been removed from the graph. In  particular cases, when in turn, $q_{v}$ or $p_{e}$ are uniform, we will specify them as node ($q$) or edge ($p$) operational probability. 
Let the degree distribution of our graph be $p_k$, i.e., a randomly chosen node has degree $k$ with probability $p_k$. We choose an edge at random and follow it to one of the nodes at its ends, then the number of edges incident on that node, other than the one we arrived
along, follows the following distribution, 
\begin{equation}
\acute{p}_k=\frac{k p_{k}}{<k>}\, ,
\end{equation}
where $<k>=\sum_{k} k p_k$,  is the average degree for the entire graph.

 Let $e_{uv}=\left\{u,v\right\}$ be the edge between node $u$ and node $v$, and $p^{*}_{uv}$ be the probability that edge $e_{uv}$  does not lead to a vertex connected via the remaining edges to the giant component. Then we have the following recursive expression,
\begin{equation}
p^{*}_{uv}=(1-p_{e_{uv}})+ p_{e_{uv}} \prod_{w} p^{*}_{vw}
\end{equation}
where product is taking over all $w$ which are not equal to $u$ and are a neighbor of $v$.
Now,  the average overall probability $p^{*}_{e}$ that a \textit{randomly} chosen edge does not belong to the giant component is,
 \begin{equation}
\label{r}
p^{*}_{e}=1-p_e+ p_e \sum_{k} \acute{p}_k (p^{*}_{e})^{k-1} :=h(p^{*}_{e})\, .
\end{equation}

The equation $p^{*}_{e}=h(p^{*}_{e})$ has the trivial solution $p^{*}_{e}=1$. Moreover if $\acute{h}(1):=\frac{d h}{d p^{*}_{e}}|_{p^{*}_{e}=1}> 1$ then $p^{*}_{e}=h(p^{*}_{e})$ has a unique root in $[0,1)$. 
This condition  $\acute{h}(1)>1$ implies that
 \begin{equation}
\label{con}
p_e > \frac{<k>}{<k^2>-<k>}\, .
\end{equation}
The existence of a giant component has a threshold at $\frac{<k>}{<k^2>-<k>}$, (the bond  percolation threshold of the graph). Note that the lower bound for $p_c$ in (\ref{con}) is meaningful, if $\frac{<k^2>}{<k>}>2$, (the Molloy-Reed criterion for existence of giant connected components,\cite{moll}).

\section{The network reliability }
\label{S:4}
Now we map the random failures process of a network  into an inverse percolation problem. The inverse percolation problem consists of finding the critical fraction of the edges for which the giant component disappears. Suppose $g$ denotes the fraction of randomly removed edges.
 The deletion of a fraction $g$ of edges corresponds to a random graph in which the
edges are occupied with probability $p = 1-g$. For small $g$, the infinite cluster which is identified as the giant component is present. The threshold for the destruction of the giant component, $g_c = 1 -p_c$,
can be thus computed from the percolation threshold at which the infinite cluster
first emerges. In this case the phase transition corresponds to the separation of a
region of damages which still allow a connected network of appreciable size from
a region in which the system is totally fragmented, \cite{past}.\\

Let $p_{e}(t)$ be the probability that edge $e$ is operational at time $t$, (the reliability of edge $e$ at time $t$). Also we define $q_{v}(t)$, to be the probability that at time $t$ the node $v \in V$ is operational. According to the previous section, if $p_{e}(t)$ is less than a critical point $p_c$, the giant cluster does not exist anymore. Define the instant at which this occurs as the lifetime of the network. 

It is known that \cite{Janson}, the binomial (uniform) random graph loses its connectivity when the number of failed edges reaches ${n \choose 2}-[M_c]$, where $M_c$ is the number of remaining edges at critical point. Since $p_c$ is a threshold for a monotone property of being connected, we have $M_c=p_c* {n \choose 2}$,\cite{Janson}.

Now using threshold $p_c$, an assessment for the network reliability at time $t$, denoted by $\hat{Rel_{c}}(G,p_e)(t)$ is given by the following equation:
\begin{equation}
\hat{Rel_{c}}(G,p_e)(t)= \sum_{i=[M_{c}]+1}^{N} \sum_{A \in [N]^{i}} \prod_{e \in A} p_{e}(t) \prod_{e \in A^{c}} [1-p_{e}(t)]  \label{r}\, ,
\end{equation}
where $N={n \choose 2}$ and $[N]^i$ stand for the family of all $i$-element subsets of $\left\{1,2,...,{n \choose 2}\right\}$.
In case when, all of $p_e(t)$'s are equal to $p(t)$, we have

\begin{equation}
\label{e2}
\hat{Rel_{c}}(G,p)(t)= \sum_{i=[M_{c}]+1}^{N} {N \choose i}  p(t)^{i} [1-p(t)]^{N-i} \, .
\end{equation}
  In \cite{daq}, authors have defined the probability that a node/edge is functional at a given time $t$ by $R(t)$. Then by the voting system based on the cite percolation critical value, they presented the network reliability at time $t$, $R_s(t)$ by the following equation:
	\begin{equation}
\label{e-old}
R_{s}(t)= \sum_{i=[n*p_{c}]+1}^{n} {n\choose i}  R(t)^{i} [1-R(t)]^{n-i} \, ,
\end{equation}
where $R(t)$ is the reliability of the generic node. In fact in their model, they assumed the same reliability for all nodes. They also gave an example in which different edges have different failure probability and then they compared the numerical results with the result from equation (\ref{e-old}). It is clear that when reliability of nodes (i.e., $R(t)$) are not same, the network reliability  cannot be obtained from equation (\ref{e-old}).\\

In the following example we will see  $\hat{Rel_{c}}(G,p)(t)$ is an $\epsilon$-approximation for the reliability polynomial \cite{ky}. In fact difference between the exact and approximate values for the reliability polynomial is zero.\\
 
\begin{example}
The reliability of complete graph $K_4$ where $p_e(t)=p(t)$ for all edges $e$ in $E$, at given time $t$ is
\begin{align}
\label{e3}
Rel(K_4,p)(t)&=p(t)^6 \quad  (F_0=1) \nonumber\\
&+ 6 p(t)^5 (1-p(t)) \quad  (F_1=6) \nonumber \\
& + 15 p(t)^4 (1-p(t))^2 \quad (F_2=15) \nonumber \\
&+ 16 p(t)^3 (1-p(t))^3 \quad (F_3=16) \, ,
\end{align}
as any sub-graph on at least $4$ edges is operational and sub-graphs with $3$ edges that are operational are those that are a spanning tree of $K_4$. There are no operational states with $2$ or fewer edges. On the other hand, for complete graph $K_4$, we have $p_c=\frac{1}{3}$. Therefore $[M_c]+1=[{4 \choose 2}*\frac{1}{3}]+1=3$. Hence we get same value for the reliability of  $K_4$ from equations (\ref{e2}) and (\ref{e3}). \\
\end{example}

The lifetime distribution of the network based on $\hat{Rel_{c}}(G,p_e)(t)$ is
\[f(t)=\frac{d(1-\hat{Rel_c}(G,p_e)(t))}{dt}\, ,\]
and the network lifetime, $T$, is given by
\[T=\int_{0}^{\infty} \hat{Rel_{c}}(G,p_e)(t) dt\, .\]
Also,
$$\hat{Rel_c}(G,p_e)(T)=p_c.$$\\

Still, the computing of $\hat{Rel_{c}}(G,p)(t)$ in (\ref{r}) by enumerating all elements in $[N]^{i}$  for large $N$, is not practicable. Thus we need to use approximation methods for computing $\hat{Rel_{c}}(G,p)(t)$. One can think of Poisson approximation method, which is used for Poisson binomial distribution( the distribution of the sum of independent and
non-identical random indicators). Fix $t>0$. Let $\mu(t)=\sum_{e \in E} p_{e}(t)$.
Then by Poisson approximation method, 
\begin{equation}
\hat{Rel_{c}}(G,p)(t) \approx \sum_{i=[M_c]+1}^{N} \frac{\mu(t)^{i} \text{e}^{-\mu(t)}}{i!}.
\end{equation}
By Le Cam's theorem \cite{le}, the approximation error for the Poisson approximation method is 
\[\sum_{k=0}^{N}| \sum_{A \in [N]^{k}} \prod_{e \in A} p_{e}(t) \prod_{e \in A^{c}} [1-p_{e}(t)]-\frac{\mu(t)^{k} \text{e}^{-\mu(t)}}{k!}| < 2\sum_{e \in E} p_{e}(t)^{2} \,.\]\\

In the following section, we look at an application of these results to some specific examples. We consider the system which can be modeled as  
\begin{itemize}
	\item [a:~] Inhomogeneous binomial random graph, and
	\item [b:~] Scale-free random graph.
\end{itemize}

\section{Inhomogeneous binomial random graph}
\label{S:5}
If edges in the binomial random graph, have different probability of being occupied, the resulting graph is called the inhomogeneous binomial random graph and denoted by $\mathcal{G}(n, \textbf{p})$, where $\textbf{p} = \left\{p_e, e\in E\right\}$, \cite{ven}. The
binomial random graph is retrieved when taking $p_e = p$ for all $e \in E$.
For this random graph, the degree distribution will be asymptotically Poisson. Hence,
$$p_k=\lambda^k \frac{\text{e}^{-\lambda}}{k!}  \, ,\quad \lambda>0\,.$$
Therefore, $p_c=\frac{<k>}{<k^2>-<k>}=\frac{1}{\lambda}=\frac{1}{<k>}\,. $\\

\begin{example}
Suppose we have a graph $G$ on $5$ nodes and the vector of edge probabilities at a given time $t$ as follows:
$$(p_{e_{1}},p_{e_{2}},p_{e_{3}},p_{e_{4}},p_{e_{5}},p_{e_{6}},p_{e_{7}},p_{e_{8}})(t)=$$
$$(\text{e}^{-0.0379*t},\text{e}^{-0.8795*t},\text{e}^{-0.7818*t},\text{e}^{-0.6949*t},\text{e}^{-0.6841*t},\text{e}^{-0.0732*t},\text{e}^{-0.1629*t},\text{e}^{-0.01045*t}).$$
In Figure~\ref{fig:gr1}, graph $G$ is presented. 
\begin{figure}[h]
	\centering
		\includegraphics [scale=0.55]{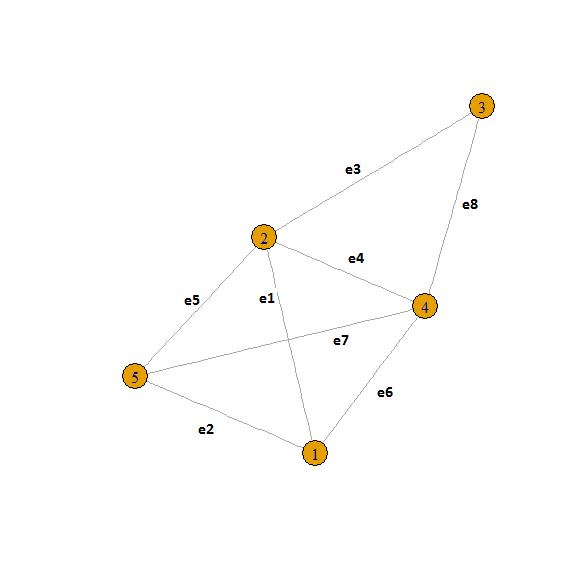}
			\caption{Graph G with 5 nodes and 8 edges.}
	\label{fig:gr1}
	\end{figure}

For this graph the value of threshold $p_c$ is $0.421053$. 
The plot of $\hat{Rel_{c}}(G,p)(t)$ verses $t, 0 \leq t \leq 15$ is shown in Figure~\ref{fig:r1}.
\begin{figure}[h]
	\centering
		\includegraphics[scale=0.52]{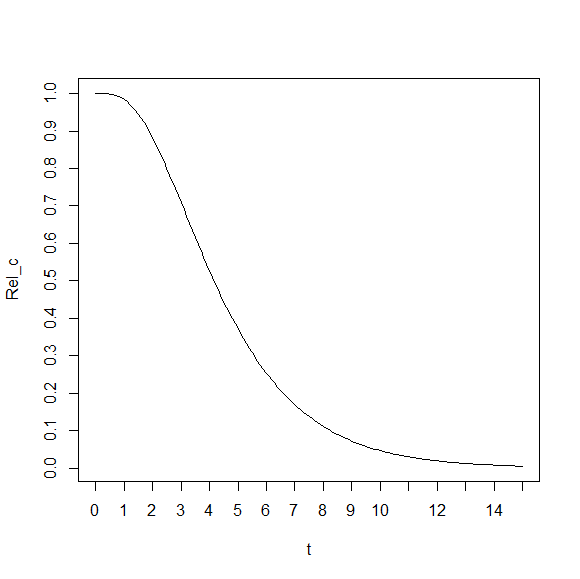}
		\caption{Plot of $\hat{Rel_{c}}(G,p)(t)$.}
		\label{fig:r1}
\end{figure}

Note that  $\hat{Rel_{c}}(G,p)(t)$ drops abruptly at $p_c=0.421053$. The corresponding $t$ for this threshold present the life time of network (between 4 and 5). In Figure~\ref{fig:gr2}, the evolution of the network $G$ is studied.
\begin{figure}[h]
	\centering
		\includegraphics [scale=0.35]{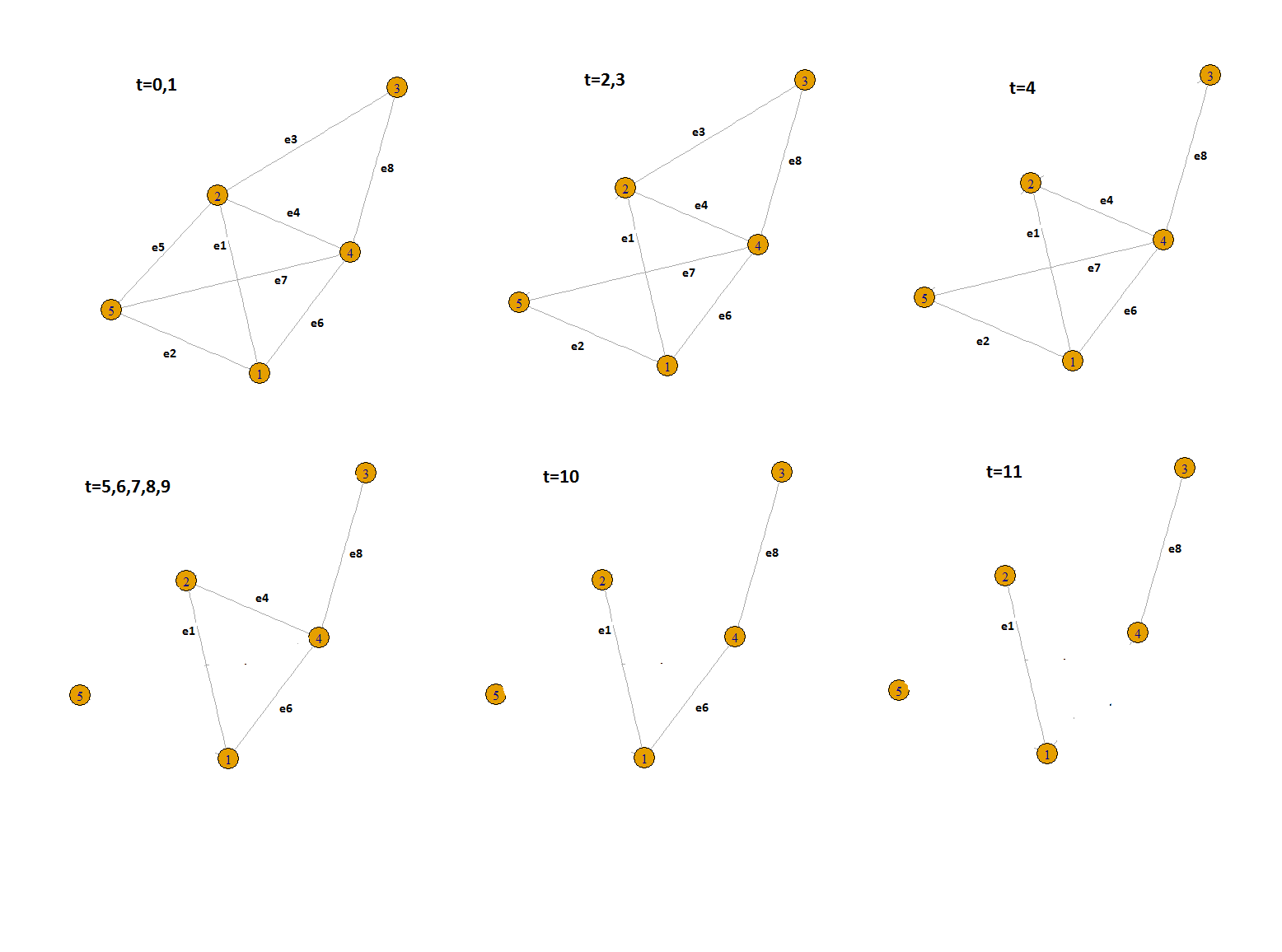}
	\caption{Evaluation of $G$ at discrete times t=0,1,...,11}
		\label{fig:gr2}
\end{figure}
\end{example}

\newpage
\section{Scale-free random graph}
\label{S:6}
In the scale-free networks majority of the nodes have few neighbors but a small portion of them have a large number of neighbors. Many of them possess a power-law degree distribution, so that the probability $p_k$ that a node has $k$ neighbors scales with $k$ as  $p_k \sim k^{-\gamma}$, where $\gamma$ is the power-law exponent \cite{bara}. Topologies of many real-world networks are thought to be scale-free, such as World Wide Web, email networks, networks of Internet routers, and protein-protein interaction networks.  
Here we consider the degree distribution as follows:

\begin{equation}
\label{power}
p_{k}= \begin{cases}
       0 & \text{ if } k=0\\
       \frac{k^{-\gamma}\text{e}^{-k/\kappa}}{\text{Li}_{\gamma}(\text{e}^{-1/\kappa})} & \text{ if } k\geq1,
      \end{cases}
\end{equation}
where $\gamma$ and $\kappa$ are constant and $\text{Li}_{\gamma}(x)$ is the polylogarithm function:
$$\text{Li}_{\gamma}(x)=\sum_{i=1}^{\infty} \frac{x^i}{i^{\gamma}}\, .$$

\begin{figure}[h]
	\centering
		\includegraphics [scale=0.7]{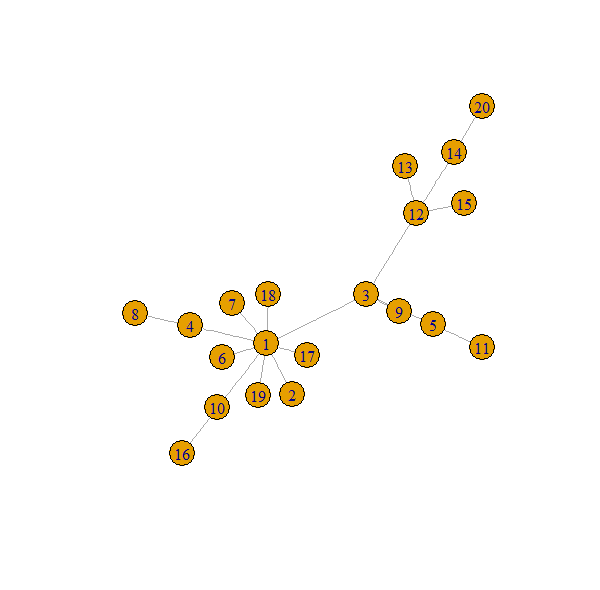}
	\caption{An example power-law graph with $\gamma=1.5$.}
	\label{fig:pow-sam}
\end{figure}

Figure~\ref{fig:pow-sam} shows a power-law graph with exponent  $\gamma=1.5$.

For power-law degree distribution, the giant component threshold is, 
\begin{equation}
p_c=\frac{<k>}{<k^2>-<k>}=\frac{\text{Li}_{\gamma-1}(\text{e}^{-1/\kappa)}}{\text{Li}_{\gamma-2}(\text{e}^{-1/\kappa})-\text{Li}_{\gamma-1}(\text{e}^{-1/\kappa})}
\label{gama}
\end{equation}

The divergence of the second moment $<k^2>$ is generally a sufficient condition to
ensure the absence of a percolation threshold on a scale-free graph. 

It is known that for $\gamma>3$ a phase transition exists at a finite $p_c$. Whereas for $2 < \gamma<3$ the transition takes place only at $p_c=0$, \cite{bara16}.\\

\begin{example}
Suppose we have a graph $G$ in which its degree distribution is as follows:
$$p_k= \frac{k^{-\gamma}}{\zeta(\gamma)}\quad ,k \geq 1 \, ,$$
where $\zeta(.)$ is the Riemann zeta function. 
Therefore from (\ref{gama}) we have
$$p_c=\frac{\zeta(\gamma-1)}{\zeta(\gamma-2)-\zeta(\gamma-1)}\,.$$
Note that  $\text{Li}_{\gamma}(1)=\zeta(\gamma).$
For many scale-free networks the degree exponent $\gamma$ is between $2$ and $3$, \cite{bara}.
The behavior of $p_c$ as a function of $\gamma$ is depicted in Figure~\ref{fig:zeta}. As the figure shows, for $\gamma \leq 3, p_c=0$ and hence the transition never takes place (unless the network is finite). Also $p_c$ is between zero and one, only in the small range $3<\gamma<3.48$ of the exponent $\gamma$.
\begin{figure}[h]
	\centering
		\includegraphics [scale=0.7]{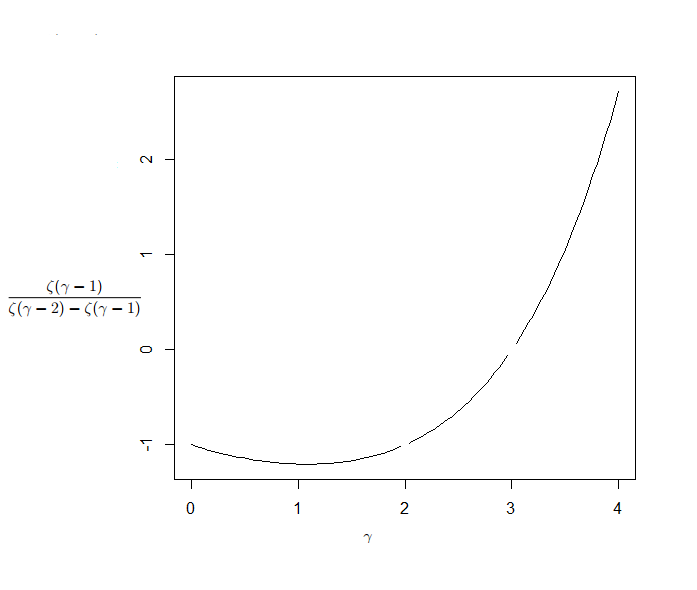}
		\caption{The function $\frac{\zeta(\gamma-1)}{\zeta(\gamma-2)-\zeta(\gamma-1)}$ for $0<\gamma <4$}
	\label{fig:zeta}
\end{figure}
As an example of real world network, we consider the internet. The degree distribution of the internet nodes follows a power law :
\begin{equation}
p_k = c k^{-\gamma}, \quad k_{\text{min}} \leq k \leq k_{\text{max}},
\end{equation}
where  $c$ is normalization constant and $k_{\text{min}}$ and $k_{\text{max}}$ are the minimal and maximal degrees, respectively, \cite{fal}. For this model we have
\begin{equation}
\label{gam}
p_{c}= \begin{cases}
       \frac{1}{\frac{\gamma-2}{3-\gamma}{k_{\text{min}}^{\gamma-2}}{k^{3-\gamma}_{\text{max}}}-1} & 2 < \gamma < 3\\
      \frac{1}{\frac{\gamma-2}{\gamma-3}(k_{\text{min}})-1} &  \gamma> 3.
      \end{cases}
\end{equation}
For the internet, $\gamma \approx 2.5$. Therefore, $p_c=\frac{1}{({k_{\text{min}}k_{\text{max}}})^{0.5}-1}\, .$  Figure~\ref{fig:pow} shows the plot of $\hat{Rel_{c}}(G,p)(t)$ when $G$ is a power-law graph with $\gamma=2.5, k_{min}=1, k_{max}=11, |E|=250 $ and edge operational probability $p_{e}(t)=\text{e}^{-0.25*t}$. As Figure~\ref{fig:pow} shows, the lifetime of the network is around $t=3.36$.\\ 
 
\begin{figure}[h]
	\centering
		\includegraphics [scale=0.7]{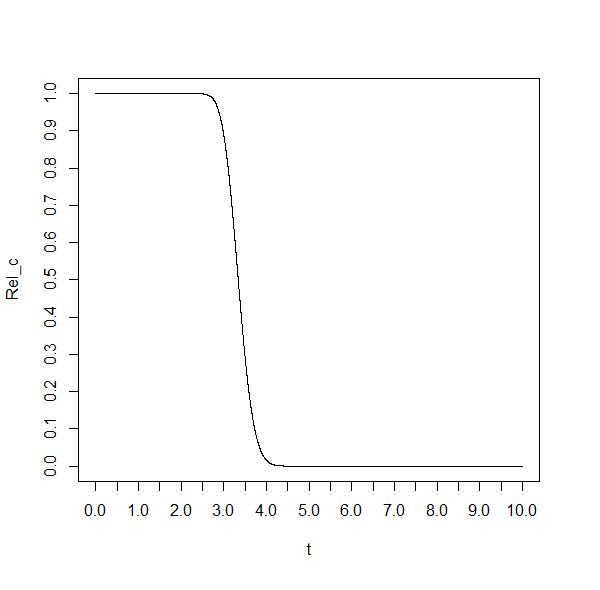}
		\caption{Plot of $\hat{Rel_{c}}(G,p)(t)$, where $G$ is power-law with $\gamma$=2.5\,.}
	\label{fig:pow}
\end{figure}

\end{example}

\bibliographystyle{amsplain}

\begin{thebibliography}{5}

\bibitem{colb1993}
Colbourn C.J., (1993). Some open problems on reliability polynomials,
,in: Congr. Numer. 93 , 187--202.

\bibitem{ch}
Larsson, C., (2014). Design of Modern Communication Networks: Methods and Applications. 1st. Academic Press.

\bibitem{colb1987}
 Colbourn, C. (1987). The Combinatorics of Network Reliability. Oxford University Press, New
York.

\bibitem{Pro1983}
 Provan J. and  Ball, M. (1983). The complexity of counting cuts and of computing the probability that a graph is connected, SIAM J. Comput. vol. 12:777--788.

\bibitem{Lomon1971}
 Lomonosov,  M. and  Polesskii,V. (1971). An upper bound for the reliability of information networks. Problems Inform. Transmission vol. 7:337--339.

\bibitem{Lomon1972}
  Lomonosov,M. and  Polesskii,V. (1972). Lower bound of network reliability, Problems Inform.
Transmission vol. 8:118--123 

\bibitem{Kru1963}
  Kruskal, J. (1963). The number of simplices in a complex, in Mathematical Optimization Techniques, R. Bellman, ed. University California Press.  251--278.

\bibitem{Stan1975}
 Stanley,  R. P. (1975). The upper bound conjecture and Cohen–Macaulay rings, Stud. Appl. Math.
vol. 14:135--142.

\bibitem{ball1982}
 Ball, M. and  Provan,J. (1982). Bounds on the reliability polynomial for shellable independence systems. SIAM J. Algebraic Discrete Methods, vol. 3:166--181.

\bibitem{ball1995}
  Ball, M.  Colbourn, C. and  Provan,J. (1995). Network reliability, in Handbook of Operations
Research: Network Models. Elsevier North-Holland),  673--762.


\bibitem{brown1996}
  Brown, J. and  Colbourn, C. (1996). Non-Stanley bounds for network reliability, J. Algebraic
Combin. vol. 5:18--36.

\bibitem{daq}
Li, D.,  Zhang, Q., Zio, E., Havlin, S. and Kang, R., (2015). Network reliability analysis based on percolation theory, Reliability. Reliability Engineering and System Safety. 142.

\bibitem{Janson}
Luczak, T. and Rucinski, A. (2000). Random Graphs., Wiley, New York. 

\bibitem{moll}
 Molloy, M. and  Reed,B. (1995).  A Critical Point for Random Graphs with a Given Degree Sequence. Random Structures and Algorithms. vol. 6:161--180.

\bibitem{past}
Pastor-Satorras, R. and Vespignani, A. (2004). Evolution and Structure of the Internet: A Statistical Physics Approach. Cambridge University Press. New York, NY, USA.

\bibitem{ky}
Kayi Lee, Hyang-Won Lee, Eytan Modiano. (2011). Reliability in Layered Networks with Random Link Failures. IEEE/ACM Transactions on Networking, vol. 19: 1835--1848.

\bibitem{le}
Le Cam, L. ,(1960), An approximation theorem for the Poisson binomial distribution. Pacific
Journal of Mathematics. vol. 10:1181--1197.

\bibitem{ven}
Remco van der Hofstad (2016). Random Graphs and Complex Networks: Volume 1. 1st. Cambridge University Press New York.

\bibitem{bara}
Albert, R. and  Barabasi, A. L. (2002). Statistical Mechanics of Complex Networks, Review Modern Physics, vol. 74: 47--97.


\bibitem{bara16}
Barabasi, A.L. and Posfai, M.(2016). Network Science. Cambridge, UK: Cambridge University Press.

\bibitem{fal}
Faloutsos, M., Faloutsos, P. and Faloutsos, C.(1999). On Power-law Relationships of the Internet Topology. ACM SIGCOMM 1999, Comput. Commun. Rev. 29, 251--262.




\end{thebibliography}

\end{document}